\documentclass[12pt,righttag]{amsart}

\usepackage{amssymb}

\newtheorem{theorem}{Theorem}
\newtheorem{remark}{Remark}
\newtheorem{lemma}{Lemma}
\newtheorem{corollary}{Corollary}
\newtheorem{conjecture}{Conjecture}
\newtheorem{proposition}{Proposition}
\newtheorem{definition}{Definition}
\newtheorem{example}{Example}

\newcommand{\astar}{\text{{\rm (}\begin{large}$\mathbf\ast$\end{large}{\rm)}}}
\newcommand{\Qed}{$\blacksquare$}

\begin{document}

\title{Isoperimetric Flow and Convexity of $H$-graphs}

\author{John McCuan}
\address{John McCuan, MSRI, 1000 Centennial Drive, Berkeley, California  94720}
\makeatletter
\email{john@msri.org}
\makeatother
\dedicatory{To Leon Simon who taught me De~Giorgi-Nash-Moser Theory, \\ and 
Craig Evans who alerted me to the fact that the time had come to use it.}
\thanks{Research supported in part by an NSF Postdoctoral Fellowship.  The 
author wishes to express thanks to the Mathematical Sciences Research 
Institute and the University of California, Berkeley Mathematics Department 
for their hospitality.}
\begin{abstract}
In this paper we consider a ``flow'' of nonparametric solutions of the 
volume constrained Plateau problem with respect to a convex planar curve.  
Existence and regularity is obtained from standard elliptic theory, and 
convexity results for small volumes are obtained as an immediate 
consequence.  Finally, the regularity is applied to show a strong stability 
condition (Theorem~\ref{stability_thm}) for all volumes considered.  
This condition, in turn, allows us to adapt an argument of Cabr\'e and 
Chanillo \cite{CabSta} which yields that any solution enclosing a non-zero 
volume has a unique nondegenerate critical point.
\end{abstract}
\maketitle

\section*{Contents}

\renewcommand{\descriptionlabel}{}

\begin{description}
\item Introduction.
\item[1.] Solvability and regularity.
\item[2.] Convexity for small volumes.
\item[3.] Stability.
\item[4.] Uniqueness of critical points.
\end{description}

\section*{Introduction}
Let $\Omega$ be a smooth ($C^{2,\alpha}$) bounded, strictly convex domain 
in ${\mathbb R}^2$.  We consider classical solutions 
$u\in C^2(\Omega)\cap C^0(\bar\Omega)$ of the constant mean curvature 
boundary value problem:
\begin{equation} 
\begin{cases}
	\text{div}\,\left(\frac{Du}{\sqrt{1+|Du|^2}}\right)  = 
		2H & \text{on $\Omega$,} \\
	u_{\big|_{\partial\Omega}}  \equiv 0 & \text{ }
\end{cases} \tag*{(\begin{large}$\mathbf\ast$\end{large})}
\end{equation}
where $H$ is a given constant.  The word ``flow'' in the title refers 
simply to the fact that we consider the solutions of 
(\begin{large}$\mathbf\ast$\end{large}) as a 
continuous family $u=u(x;H)$ for $x\in \bar\Omega$ and $H$ in a suitable 
interval (which we can think of as a time interval).  We remark also that 
each solution (for fixed $H$) provides the least area graph which, along 
with the domain $\Omega$, encloses a certain prescribed volume.  See \S3 
below for details.

In the following \S1 we delineate existence and regularity results for the 
family $u(x;H)$.  These results follow from standard theorems in elliptic 
pde for which we cite extensively \cite{GilEll}.  The utility of the results 
in \S1 is demonstrated throughout the rest of the paper.

In \S2 we characterize (in terms of Poisson's equation) two aspects of the 
convexity properties of solutions that enclose a sufficiently small volume.  
It follows from this characterization, for example, that if $\Omega$ has 
boundary an ellipse, then for $|H|\ne 0$ small enough, solutions of 
(\begin{large}$\mathbf\ast$\end{large}) are convex.  It is also shown for 
all ($C^{2,\alpha}$) convex domains that 
small enough solutions (in the sense just mentioned) have convex level 
curves.  More precisely, we show that small solutions are 
$1/2$ (power) convex (see \S2 for details).

We are unable at present to derive any convexity property for larger 
solutions.  In \S3, however, we establish a strong stability 
condition.  The condition is ``strong'' in the sense that it holds 
for arbitrary variations, not just those that preserve the volume 
with respect to which the solution is known to be a minimizer of area.  
The stability condition is then used in \S4 to show that each non-zero 
solution has a unique critical point.  In this regard, we follow an 
argument of Cabr\'e and Chanillo \cite{CabSta} which they applied to 
semilinear equations.  In particular, each solution has a unique global 
extremum (maximum or minimum), a corollary that can also be deduced from 
the convexity of the level curves for small solutions (\S2).  It also 
follows from the method of Cabr\'e and Chanillo that the critical point 
is nondegenerate, i.e., solutions are strictly convex in some neighborhood 
of the critical point.

Most of the arguments below (usually in an ill-formed state) and many 
ill-fated versions of them have been inflicted on my friends and 
colleagues.  I thank them for their patience and helpful comments.  Among 
them are Claire Chan, Mikhail Feldman, Melinda McCuan, Robert Osserman, 
and Tatiana Toro.  I should specifically like to thank Robert Finn 
for Example~\ref{nodex}, Henry Wente for suggesting the ``form'' of 
Theorem~\ref{1convex}, and Brian White for making Remark~\ref{velocity}.  
Finally, I am indebted to 
David Hoffman for introducing me to ``bubble problems.''

\section{Solvability and regularity}\label{s1}

The results of this section apply more generally to problems 
\begin{equation} 
\begin{cases}
	\text{div}\,\left(\frac{Du}{\sqrt{1+|Du|^2}}\right)  = 
		2H & \text{on $\Omega$,} \\
	u_{\big|_{\partial\Omega}}  \equiv \phi & \text{ }
\end{cases} \tag*{(\begin{large}$\mathbf\ast\ast$\end{large})}
\end{equation}
where $\Omega$ is a $C^{2,\alpha}$ domain in ${\mathbb R}^n$ 
(not necessarily convex) and either (i) $\phi \equiv 0$, or (ii) 
$\phi$ (a function defined on $\partial\Omega$) extends to a $C^{2,\alpha}$ 
function on $\bar\Omega$ and the mean curvature of $\partial\Omega$ is 
everywhere positive.  The proofs, for the most part, apply to both cases 
(i) and (ii) though certain details not directly related to our main 
results are presented in Appendix~\ref{appthm1}.

We begin by determining a suitable interval on which to consider the mean 
curvature $H$.

\begin{theorem}[solvability]\label{solvabilitythm}  There is a unique value 
$H_{\text{max}}>0$ depending on $\Omega$ such that the following hold.
\begin{description}
\item[{\rm i.}] If $|H| \le H_{\text{max}}$, {\em \astar} has a unique solution $u=u(x;H)$.
\item[{\rm ii.}] If $|H| > H_{\text{max}}$, {\em \astar} has no solution.
\end{description}
Furthermore, if $|H| < H_{\text{max}}$, $u\in C^{2,\alpha}(\bar\Omega)$, but 
\begin{equation*}
	\sup_{x\in\Omega} |Du(x;H_{\text{max}})| = \infty.
\end{equation*}
\end{theorem}

\begin{remark} \label{rem1}  
If $\kappa$ denotes the minimum mean curvature of $\partial\Omega$ 
and $2|H| \le \kappa$, then solvability follows from Theorem~16.11 
\cite[pg. 409]{GilEll}.  For the particular boundary condition in {\em \astar} 
solvability will, in general, persist for $2|H| > \kappa$, and we wish 
to include these solutions in our discussion.
\end{remark}
\begin{remark}
Any results numbered 6.2-16.11 without specific reference, refer to 
\cite{GilEll}.
\end{remark}
\begin{remark}\label{rem2}
The mean curvature operator that appears in {\em \astar} will be denoted by 
${\mathcal M}$.  Furthermore, given $p\in {\mathbb R}^n$ we define the 
vector $A=p/\sqrt{1+|p|^2}$ and write ${\mathcal M}$ in its 
``pure divergence form'' {\em (\ref{pdform})} which expands to an alternative 
``non-divergence quasilinear form'' {\em (\ref{ndql})}.
\begin{equation}
\begin{align}
	{\mathcal M}u &= \sum_i D_i A^i(Du) \label{pdform} \\
		  &= \sum_{i,j} \frac{\partial A^i}{\partial p_j}(Du) 
			D_iD_ju.\label{ndql}
\end{align}
\end{equation}
The coefficients in the last expression will be denoted by 
$A_{ij} = A_{ij}(Du)$, and one easily checks that the coefficient matrix 
$(A_{ij})$ is positive definite, i.e., the operator is 
elliptic (and uniformly elliptic if $|Du|$ remains bounded).  

We may also consider ${\mathcal M}$ in ``non-divergence linear form'' 
\begin{equation}\label{linform}
	{\mathcal M}u = \sum_{i,j} a_{ij} D_iD_ju,
\end{equation}
by simply setting $a_{ij}(x) = A_{ij}(Du(x))$.  

Given a solution $u\in C^2(\bar\Omega)$ of {\em \astar}, it follows that 
$u\in C^\omega(\Omega)\cap C^{2,\alpha}(\bar\Omega)$.  See Corollary~16.7 
\cite[pg. 407; see also pp. 109--111]{GilEll}.
\end{remark}
\begin{remark}\label{rem3}
When convenient, we indicate the dependence of the problem {\em \astar} on 
the mean curvature by a subscript: $\astar_H$.
\end{remark}

\noindent{\bf Proof of Theorem~\ref{solvabilitythm}.}  We set 
\begin{equation*}
	H_{\text{max}} \equiv \sup\{ H: {\text{\astar \ is (uniquely) 
		solvable in $C^2(\bar\Omega)$}} \}.
\end{equation*}
Let us assume for the moment that $H_{\max} > 0$.  Note that in our 
case of primary interest (convex planar domains) this condition of 
``nondegeneracy'' follows from Remark~\ref{rem1}.  Uniqueness of a given 
solution $u$ is immediate from the comparison principle, Theorem~10.1 
\cite[pg. 263]{GilEll}; see also Theorem~10.2.

Following the usual Leray-Schauder approach 
(Theorem~13.8 \cite[pg. 331]{GilEll}), solvability follows if there is a 
constant $M$ such that the apriori bound
\begin{equation}\label{b0}
	|u|_{C^1(\Omega)} = \sup_\Omega |u| + \sup_\Omega |Du| < M
\end{equation}
holds for any $C^{2,\alpha}(\bar\Omega)$ solution of $\astar_{\sigma H}$.  
Note that the constant $M$ is required to be independent of $u$ and 
$\sigma>0$.

To obtain a such a bound we assume that $|H| < H_{\max}$ 
and take $\tilde{H} \in (|H|, H_{\max})$ such that $\astar_{\tilde{H}}$ 
has a solution $\tilde{u}\in C^2(\bar\Omega)$.  ($\tilde{H}$ exists by 
the definition of $H_{\max}$.)  By the 
comparison principle, any solution of $\astar_{\sigma H}$ satisfies 
$|u| \le |\tilde{u}|$ and thus, 
\begin{equation}\label{b1}
	|Du|_{\big|_{\partial\Omega}} \le 
			|D\tilde{u}|_{\big|_{\partial\Omega}} \le 
			|\tilde{u}|_{C^1(\Omega)} 
					\equiv \tilde{M} < \infty.
\end{equation}

On the other hand, we can differentiate the expression 
(\ref{linform}) to obtain a linear elliptic equation satisfied by 
$v = D_ku$.  In fact, 
\begin{equation}\label{linearize1}
	{\mathcal L}v \equiv \sum_{i,j}a_{ij} D_iD_jv + \sum_l 
		b_l D_lv = 0 
\end{equation}
where
\begin{equation*}
	b_l(x) = \sum_{i,j} 
		\frac{\partial^2A^i}{\partial p_j\partial p_l}
			(Du) D_iD_ju.
\end{equation*}
By the weak maximum principle, Theorem~3.1 \cite[pg. 32]{GilEll}, 
\begin{equation*}
	\sup_\Omega|D_ku| = \sup_{\partial\Omega}|D_ku|.
\end{equation*}
Consequently, we have from (\ref{b1})
\begin{equation*}
\begin{align}
	\sup_\Omega|Du| &\le \sqrt{n} \max_k \sup_\Omega|D_ku| \notag \\
			&=   \sqrt{n} \max_k \sup_{\partial\Omega}|D_ku| 
					\notag \\
			&\le \sqrt{n} \tilde{M}. 
\end{align}
\end{equation*}

We have therefore established the apriori bound (\ref{b0}) with 
$M = (1 + \sqrt{n}) \tilde{M}$, and solvability 
follows for $|H| < H_{\max}$.  

Since our primary results concern convex planar domains and the case 
$|H| < H_{\max}$, we postpone the remainder of the proof of 
Theorem~\ref{solvabilitythm} (see Appendix~\ref{appthm1}) and use only the 
assertions established above.  Concerning the extremal solution 
$u(x;H_{\max})$, we note that the nondegeneracy condition, $H_{\max} > 0$, 
and the gradient blow-up condition follow in general from a ``short time 
existence'' result, Theorem~\ref{shorttime}, that is of independent 
interest.  $\square$

\medskip

This is a convenient time to point out two other immediate consequences 
of the comparison principle (Theorem~10.1).
\begin{corollary}[monotonicity and symmetry]\label{monotonicity}
If $-H_{\max} < H \le \tilde{H} < H_{\max}$, then 
\begin{equation*}
	u(x;H) \ge u(x,\tilde{H})
\end{equation*}
for all $x\in\Omega$ with equality only if $H=\tilde{H}$.  There also 
holds
\begin{equation*}
u(x,-H) \equiv -u(x,H).
\end{equation*}
\end{corollary}

\medskip

The uniform estimate (\ref{b0}) in the proof of Theorem~\ref{solvabilitythm} 
 allows us to 
concentrate on certain questions of uniformity in $H$ for higher 
derivatives of $u$ and to effectively ignore dependencies on 
ellipticity constants and bounds for the top order coefficients 
(usually denoted by $\lambda$ and $\Lambda$ respectively in \cite{GilEll}).  
This observation is recorded for reference in the following
\begin{corollary}\label{ellipcor}
If $0 < \tilde{H} < H_{\max}$, then we have a uniform bound 
\begin{equation}\label{b0cor}
	|u|_{C^{1,\tilde{\alpha}}(\Omega)} \le \tilde{M}
\end{equation}
depending only on $\tilde{H}$.  In addition, the coefficients in 
{\em (\ref{ndql})} and {\em (\ref{linform})} are uniformly 
elliptic and bounded for all $x\in \bar\Omega$ and $|H|<\tilde{H}$, 
i.e., there is some $\tilde{M}$ and some $\tilde{\lambda} > 0$, 
both independent of $H$, such that
\begin{equation}\label{coefb}
	|a_{ij}|_{C^{0,\tilde{\alpha}}(\Omega)} \le \tilde{\Lambda},
		\ \text{and} 
\end{equation}
\begin{equation}\label{11}
	\sum_{i,j} a_{ij}(x) \xi_i \xi_j \ge \tilde{\lambda} |\xi|^2
\end{equation}
for all $\xi \in {\mathbb R}^n$ and $x\in\bar\Omega$.
\end{corollary}

\noindent{\bf Proof.}  The ellipticity constant $\lambda$ for $(a_{ij})$ is 
$\inf_\Omega(1+|Du|^2)^{-3/2}$.  Thus, by (\ref{b0}) $\lambda \ge 
(1+n\tilde{M}^2)^{-3/2} \equiv \tilde{\lambda} >0$, and (\ref{11}) holds.  

Since $a_{ij} = A_{ij}(Du)$ and $A$ is smooth, $|a_{ij}|_{C^0(\Omega)} \le 
|A|_{C^2(G)}$ where $G = \{ Du \in {\mathbb R}^n : x\in \Omega \}$.  Thus, 
we have a bound, $|a_{ij}|_{C^0(\Omega)} \le \tilde{\Lambda}_0$.

Given $\tilde{\lambda}$, $\tilde{\Lambda}_0$ and $M$ from (\ref{b0}), 
Theorem~13.2 \cite[pg. 323]{GilEll} implies\footnote{Strictly speaking, 
	Gilbarg and Trudinger state the theorem with a dependence on 
	$K\equiv |u|_{C^1(\Omega)}$.  The dependence however essentially 
	arises from the De~Giorgi-Nash estimates (Theorems~8.22 and 8.29 
	\cite[pp. 200--205]{GilEll}) when they are applied to the equation 
	\begin{equation*}
		\sum_{ij}D_i(a_{ij} D_jw) = 0
	\end{equation*}
	where $w=D_ku$.  In this instance, only a bound for 
	$|a_{ij}|_{C^0(\Omega)}$ is required, not the explicit value $K$.}  
that for some $\tilde{\alpha} = \tilde{\alpha}(\tilde{\lambda}, 
\tilde{\Lambda}_0,\Omega)$
\begin{equation*}
		[Du]_{C^{\tilde{\alpha}}(\Omega)} = 
			\sup_{x,y\in\Omega, x\ne y} 
			\frac{|Du(x)-Du(y)|}{|x-y|} \le C
\end{equation*}
where $C=C(\tilde{\lambda},\tilde{\Lambda}_0,M,\Omega,\phi)$.  Thus, 
(\ref{b0cor}) holds.  Notice that we returned to the ``pure divergence 
form'' ${\mathcal M}u = \text{div}\, A(Du) = 2H$ in order to apply 
Theorem~13.2.  

Extending slightly our estimate for $a_{ij}$, see Lemma~\ref{A1} 
Appendix~\ref{apphol}, we see that $[a_{ij}]_{C^{\tilde{\alpha}}(\Omega)} 
\le |A_{ij}|_{C^1(G)} [Du]_{C^{\tilde{\alpha}}(\Omega)} \le |A|_{C^3(G)} 
\tilde{M}$.  Thus, (\ref{coefb}) follows from (\ref{b0cor}).  \Qed

\smallskip

\begin{remark} The regularity assertion of Remark~\ref{rem2}, the 
solvability theorem above (n.b. Theorem~13.8), and the bounds given 
in Corollary~\ref{ellipcor} all depend crucially on the 
$C^\alpha$ gradient bound of Ladyzhenskaya and Ural'tseva 
(Theorem~13.2) which follows from the 
De~Giorgi-Nash-Moser Theory of Chapter~8 (or alternatively---in two 
dimensions---from earlier 
results of Morrey).  This dependence also presents 
itself as the main difficulty when one tries to prove the continuity theorem 
below by applying the Schauder estimates from Chapters~6 and 8. 
\end{remark}

\begin{theorem}[continuity]\label{ctythm}
If $|H_0|,|H| \le \tilde{H} < H_{\max}$, then the following estimates hold.
\begin{description}
\item[{\rm (i)}] If $\Omega$ is a $C^{k,\alpha}$ domain for some $k \ge 1$ 
and $\alpha>0$, 
\begin{equation*}
	|u(\,\cdot\,;H) - u(\,\cdot\,;H_0)|_{C^{k,\alpha_0}(\Omega)} 
			\le C_0 |H-H_0|
\end{equation*}
where $C_0$ and $\alpha_0 > 0$ are independent of $H$.
\item[{\rm (ii)}] If $\Omega'\subset\subset\Omega$, then 
\begin{equation*}
	|u(\,\cdot\,;H) - u(\,\cdot\,;H_0)|_{C^k(\Omega')} 
			\le C'_0 |H-H_0|
\end{equation*}
where $C'_0$ depends on $k$ and $\Omega'$ but not $H$.
\end{description}
\end{theorem}

\noindent{\bf Proof.}  Let us consider $H_0$ and $\tilde{H} < H_{\max}$ 
fixed and write $u_0 = u(\,\cdot\,,H_0)$, so that our primary focus 
becomes dependence on $H$.  Accordingly, 
we will denote various constants that are independent of $H$ with a 
subscript $0$, and by $C_0$ in particular.  As mentioned above, constants 
that can be taken to depend only on $\tilde{H}$ will be, for the most part, 
ignored.  Various other constants will be denoted by $C$.

We begin by observing that the difference 
$v = u-u_0 = u(\,\cdot\,;H) - u(\,\cdot\,;H_0)$ 
satisfies a linear elliptic equation.  
If $\delta = H-H_0$, 
\begin{align*}
2\delta &= \text{div}\, A(Du)-\text{div}\, A(Du_0) \\
	&= \sum_i D_i\left\{ A^i(tDu + (1-t)Du_0)_{\big|_{t=0}^1} \right\} \\
	&= \sum_i D_i(\alpha_{ij}(x)\, D_jv)
\end{align*}
where 
\begin{equation}\label{alphij}
	\alpha_{ij} = \alpha_{ij}^\delta(x) = 
		\int_0^1 \frac{\partial A^i}{\partial p_j}
					(tDu + (1-t)Du_0)\, dt.
\end{equation}

The uniform ellipticity of Corollary~\ref{ellipcor} is easily seen to 
hold for $(\alpha_{ij})$.  Thus, we see that $v$ satisfies the uniformly 
elliptic divergence structure boundary value problem
\begin{equation} 
\begin{cases}\label{diffeqn}
	\sum_{i,j}D_i(\alpha_{ij}D_jv) = 
		2\delta & \text{on $\Omega$,} \\
	v_{\big|_{\partial\Omega}}  \equiv 0, & \text{ }
\end{cases} 
\end{equation}
and the estimates of the theorem follow, at least formally, from slight 
extensions of the ``weak'' Schauder estimates Theorems~8.33 and 8.32 
\cite[pg. 210]{GilEll}.  For reference we give the statements as they 
apply to a general divergence form linear boundary value problem.

\smallskip

\begin{quote}{\em 
Let $v\in C^{k,\alpha}(\Omega)$ (where $\Omega$ is a bounded domain and 
$k = 1$ or $2$ or $3\ldots$) be a (weak) solution of the linear boundary 
value problem
\begin{equation*} 
\begin{cases}
	Lv = g + \sum_jD_jf_j & \text{on $\Omega$,} \\
	v_{\big|_{\partial\Omega}}  \equiv \phi, & \text{ }
\end{cases} 
\end{equation*}
where $Lv = \sum_{i,j}D_i(\alpha_{ij}(x) D_jv + \beta_i(x) v) + 
\sum_j c_j(x) D_jv +d(x) v$ and the coefficients satisfy 
\begin{equation*}
	|\alpha_{ij}|_{C^{k-1,\alpha}(\Omega)},\, 
	|\beta_i|_{C^{k-1,\alpha}(\Omega)},\, 
	|c_j|_{C^{k-1}(\Omega)},\,
	|d|_{C^{k-1}(\Omega)} \le \Lambda_k.
\end{equation*}
{\bf Theorem $\mathbf{\text{\bf 8.32}}'$.}  If 
$\Omega' \subset\subset \Omega$, then 
\begin{equation*}
	|v|_{C^{k,\alpha}(\Omega')} \le C(|v|_{C^0(\Omega)} 
		+ |g|_{C^{k-1}(\Omega)} 
		+ |f|_{C^{k-1,\alpha}(\Omega)})
\end{equation*}
where $C = C(n,\lambda,\Lambda_k,\Omega',\Omega)$ ($\lambda$ being the 
ellipticity constant for $(\alpha_{ij})$) and $|f|_{C^{k-1,\alpha}(\Omega)} 
= \sum |f_j|_{C^{k-1,\alpha}(\Omega)}$.
\begin{equation*}
	\text{    }
\end{equation*}
{\bf Theorem $\mathbf{\text{\bf 8.33}}'$.}  If 
$\Omega$ is a $C^{k,\alpha}$ domain and $v\in C^{k,\alpha}(\bar{\Omega})$, 
then 
\begin{equation*}
	|v|_{C^{k,\alpha}(\Omega)} \le C(|v|_{C^0(\Omega)} 
		+ |\phi|_{C^{k,\alpha}(\Omega)}
		+ |g|_{C^{k-1}(\Omega)} 
		+ |f|_{C^{k-1,\alpha}(\Omega)})
\end{equation*}
where $C = C(n,\lambda,\Lambda_k,\Omega)$.
}
\end{quote}

\smallskip

\noindent Applying Theorem~${\text{8.33}'}$, for example, 
we have:  {\em If $|\alpha_{ij}|_{C^{k-1,\alpha}(\Omega)} \le \Lambda$, then 
\begin{equation*} 
	|v|_{C^{k,\alpha}(\Omega)} \le C( |v|_{C^0(\Omega)} + 2|\delta|)
\end{equation*}
where $C = C(n,\lambda,\Lambda,\Omega)$.}

The only dependence in $C$ on $H$ is through $\Lambda$, so we see 
that the following two lemmas together establish statement (i) of the 
theorem.
\begin{lemma}\label{cklem}
There are constants $\Lambda_0$ and $\alpha_0$ (independent of $H$) 
such that 
\begin{equation*} 
	|\alpha_{ij}|_{C^{k-1,\alpha_0}(\Omega)} \le \Lambda_0.
\end{equation*}
\end{lemma}
\begin{lemma}\label{zerolem}
There is some $C_0>0$ 
such that 
\begin{equation*} 
	|u(\,\cdot\,;H) - u(\,\cdot\,;H_0)|_{C^0(\Omega)} 
			\le C_0 |\delta|.
\end{equation*}
\end{lemma}
Lemma~\ref{zerolem} follows immediately from Theorem~8.16 
\cite[pg. 191]{GilEll}.  In fact, we have 
\begin{equation*}
	|v|_{C^0(\Omega)} \le (C/\lambda)(2|\delta|) |\Omega|^{2/q}
\end{equation*}
for any $q > n$ where $C = C_0(n,q,\Omega)$ is independent of $H$.

\noindent{\bf Proof of Lemma~\ref{cklem}.}  We see from the definition 
of $\alpha_{ij}$ in (\ref{alphij}) that for any $\alpha >0$, 
\begin{equation*}
	|\alpha_{ij}|_{C^{k-1,\alpha}(\Omega)} \le \sup_{0\le t\le 1} 
			|\tilde{\alpha}_{ij}|_{C^{k-1,\alpha}(\Omega)}
\end{equation*}
where $\tilde{\alpha}_{ij} = \tilde{\alpha}_{ij}(x,t) \equiv 
A_{ij}(tDu + (1-t)Du_0)$.  By Lemma~\ref{A1} (see Appendix~\ref{apphol}) 
$|\tilde{\alpha}_{ij}|_{C^{k-1,\alpha}(\Omega)}$ can be bounded in terms 
of $B_1$ and $B_2$ where 
\begin{equation*} 
	|A_{ij}|_{C^k(G)} \le B_1,\ \text{and}
\end{equation*}
\begin{equation*} 
	|tDu + (1-t)Du_0|_{C^{k-1,\alpha}(\Omega)} \le B_2.
\end{equation*}
In this instance 
$G = \{ tDu + (1-t)Du_0 \in {\mathbb R}^n : x\in\Omega \}$.  
As in the estimates in Corollary~\ref{ellipcor}, $A_{ij}$ is smooth and 
$G$ is bounded independently of $H$, so $B_1$ can be taken independently 
of $H$.

To find $B_2$ independently of $H$ it suffices to bound 
$|u|_{C^{k,\alpha}(\Omega)}$.  We proceed by induction.

The initial case $k=1$ is obtained from (\ref{b0cor}) by taking 
$\alpha_0 = \tilde{\alpha}$.

For $k\ge 2$, we take $\alpha_0$ to be the minimum of 
$\tilde{\alpha}$ and the H\"older exponent of $\partial\Omega$ and 
assume inductively that 
\begin{equation*}
	|u|_{C^{k-1,\alpha_0}(\Omega)}  \le M_0,
\end{equation*}
and (as a consequence)
\begin{equation*}
	|\alpha_{ij}|_{C^{k-2,\alpha_0}(\Omega)} \le \Lambda_0
\end{equation*}
for some $\Lambda_0$ independent of $H$.  The latter assumption puts 
us in a position to apply an extension of the ``classical'' Schauder 
global estimate, Theorem~6.6 n.b., Problem~6.2, which again we state 
for convenience.

\smallskip

\begin{quote}{\em 
Let $v\in C^{k,\alpha}(\Omega)$ be a (classical) solution of the 
linear boundary value problem
\begin{equation*} 
\begin{cases}
	Lv = f & \text{on $\Omega$,} \\
	v_{\big|_{\partial\Omega}}  \equiv \phi, & \text{ }
\end{cases} 
\end{equation*}
where $Lv = \sum_{i,j} a_{ij}(x) D_iD_jv + \sum_i b_i(x) D_iv + c(x) v$ 
and the coefficients satisfy 
\begin{equation*}
	|a_{ij}|_{C^{k-2,\alpha}(\Omega)},\, 
	|b_i|_{C^{k-2,\alpha}(\Omega)},\, 
	|c|_{C^{k-2,\alpha}(\Omega)} \le \Lambda_k.
\end{equation*}
{\bf Theorem $\mathbf{\text{\bf 6.6}}'$.}  If 
$\Omega$ is a $C^{k,\alpha}$ domain and $v\in C^{k,\alpha}(\bar{\Omega})$, 
then 
\begin{equation*}
	|v|_{C^{k,\alpha}(\Omega)} \le C(|v|_{C^0(\Omega)} 
		+ |\phi|_{C^{k,\alpha}(\Omega)}
		+ |f|_{C^{k-2,\alpha}(\Omega)})
\end{equation*}
where $C = C(n,\lambda,\Lambda_k,\Omega)$.
}
\end{quote}

\smallskip

\noindent When applied to the equation in $\astar$ Theorem~${\text{6.6}}'$ 
yields
\begin{align*}
	|u|_{C^{k,\alpha_0}(\Omega)} &\le 
		C_0 (|u|_{C^0(\Omega)} + 2|H|) \\
	&\le C_0 \quad \text{independent of $H$.  }
\end{align*}

The induction is concluded with the use of Lemma~\ref{A1} which implies 
a bound for $|\alpha_{ij}|_{C^{k-1,\alpha_0}(\Omega)}$.  This establishes 
Lemma~\ref{cklem} and Theorem~\ref{ctythm} part (i).

\medskip

If we replace $\Omega$ by $\Omega'\subset\subset\Omega$ in the proof of 
Lemma~\ref{cklem} and use 

\smallskip

\begin{quote}{\em 
{\bf Theorem $\mathbf{\text{\bf 6.2}}'$.}  Let $v\in C^{k,\alpha}(\Omega)$ 
be a solution of the linear boundary value problem described just above.  
Then
\begin{equation*}
	|v|_{C^{k,\alpha}(\Omega')} \le C(|v|_{C^0(\Omega)} 
		+ |f|_{C^{k-2,\alpha}(\Omega)})
\end{equation*}
where $C = C(n,\lambda,\Lambda_k,\Omega',\Omega)$.
}
\end{quote}

\smallskip

\noindent then the same reasoning yields an estimate
\begin{equation}\label{nicest}
	|u|_{C^{k,\alpha_0}(\Omega')} \le C_0(k),
\end{equation}
which in turn gives by Lemma~\ref{A1}

\medskip

\noindent{\bf Lemma $\mathbf{\text{\bf 1}}'$.}  For any $k$, there is 
a constant $\Lambda'_0$ such that 
\begin{equation*}
	|\alpha_{ij}|_{C^k(\Omega')} \le \Lambda'_0.
\end{equation*}
Having made this observation, Theorem~\ref{ctythm} part (ii) follows 
from Theorem~${\text{8.32}'}$.  \Qed

\medskip

The main theorem of this section is the following:
\begin{theorem}[regularity]\label{regularitythm}
$u\in C^\infty(\Omega \times (-H_{\max},H_{\max}))$.
\end{theorem}
\noindent{\bf Proof.}  By Theorem~\ref{ctythm} statement (i), if 
$\tilde{H} < H_{\max}$ and $|H| + |\delta| \le \tilde{H}$, then 
\begin{equation*}
	\Delta_H^\delta u \equiv 
	\frac{u(\,\cdot\,;H+\delta)-u(\,\cdot\,;H)}{\delta} = 
		\frac{u_1-u_0}{\delta}
\end{equation*}
satisfies
\begin{equation*}
	|\Delta_H^\delta u|_{C^{2,\alpha}(\bar{\Omega})} \le C_0,
\end{equation*}
i.e., $\{ \Delta_H^\delta u \}_{|H|+|\delta|\le\tilde{H}}$ is bounded in 
$C^{2,\alpha}(\bar{\Omega})$.  By Lemma~6.36, this set is therefore 
precompact in $C^2(\bar{\Omega})$.  It follows that there is a function 
$\dot{u} \in C^2(\bar{\Omega})$ such that $\lim_{\delta\to 0} 
|\dot{u} - \Delta^\delta_H u|_{C^2(\bar{\Omega})} = 0$, and 
\begin{equation}
\begin{cases}\label{dotu}
	\sum_{ij} D_i \left( a_{ij} D_j\dot{u} \right)  = 
		2& \text{on $\Omega$,} \\
	\dot{u}_{\big|_{\partial\Omega}}  \equiv 0.& \text{ }
\end{cases}
\end{equation}

Technically, $\Delta^\delta_H u$ satisfies a boundary value problem 
(divide equation (\ref{diffeqn}) by $\delta$), and by taking the limit 
of subsequences $\delta_j \to 0$ we arrive at (\ref{dotu}).  The 
existence of the limit as $\delta \to 0$ then follows from the 
uniqueness of solutions to (\ref{dotu}).

It also follows from (\ref{dotu}) that 
$\dot{u} \in C^\infty(\Omega) \cap C^{k,\alpha}(\bar{\Omega})$ (if 
$\Omega$ is a $C^{k,\alpha}$ domain).  
\begin{remark}\label{velocity}
Note that the velocity $\dot{u}$ of our flow is determined by a linear 
boundary value problem---rather than by a local expression as for 
example in the heat equation.
\end{remark}

In our derivation of (\ref{dotu}) we assumed that $\Omega$ was at 
least $C^{2,\alpha}$.  So as not to require $\partial\Omega$ to be 
inordinately smooth, we work locally from now on.

We next consider the continuity of $\dot{u}$ as a function of $H$.  
Extending the notation above, we write $\dot{v} = \dot{u} - \dot{u}_0 \equiv 
\dot{u}(\,\cdot\,;H) - \dot{u}(\,\cdot\,;H_0)$.  Note that $\dot{v}$ 
satisfies a linear boundary value problem
\begin{equation}\label{dotv}
\begin{cases}
	\sum_{ij} D_i \left( a_{ij}^0 D_j\dot{v} \right)  = 
		f & \text{on $\Omega$,} \\
	\dot{v}_{\big|_{\partial\Omega}}  \equiv 0& \text{ }
\end{cases}
\end{equation}
where $a_{ij}^0 = A_{ij}(Du_0)$ and 
$f = - \sum_{i,j} D_i[(a_{ij} - a_{ij}^0) D_j\dot{u}]$.  

Letting $\Omega'$ be a smooth domain compactly contained in $\Omega$ 
(n.b. Problem 6.9) and $\Omega'' \subset\subset \Omega'$, we have from 
Theorem~$\text{8.32}'$ 
\begin{equation*}
	|\dot{v}|_{C^{k,\alpha}(\Omega'')} \le C(|\dot{v}|_{C^0(\Omega')} 
		+ \sum_{i,j}|f_{ij}|_{C^{k-1,\alpha}(\Omega')})
\end{equation*}
where $C = C(n,\lambda,\Lambda,\Omega'',\Omega')$ and $f_{ij} = 
(a_{ij} - a_{ij}^0) D_j\dot{u}$.  In this case, $C$ is independent of 
$H$.  Thus, we have an estimate 
\begin{equation}\label{C}
	|\dot{v}|_{C^k(\Omega'')} \le C''_0|H-H_0| = C''_0|\delta|
\end{equation}
(for arbitrary $k$ as in Theorem~\ref{ctythm})
as long as the following lemmas hold.
\begin{lemma}\label{A}
	There is some $C_0$ and some $\alpha_0>0$ such that
\begin{equation*}
	|f_{ij}|_{C^{k-1,\alpha_0}(\Omega')} \le C_0|\delta|.
\end{equation*}
\end{lemma}
\begin{lemma}\label{B}
	There is some $C_0$ such that
\begin{equation*}
	|\dot{v}|_{C^0(\Omega')} \le C_0|\delta|.
\end{equation*}
\end{lemma}

As before Theorem~8.16 implies 
\begin{equation*}
	|\dot{v}|_{C^0(\Omega')} \le C(q)\|f\|_{L^{q/2}}
\end{equation*}
for any $q>n$.  Therefore, Lemma~\ref{B} follows from Lemma~\ref{A}.

\smallskip

\noindent{\bf Proof of Lemma~\ref{A}.}  Since 
$|a_{ij} - a_{ij}^0|_{C^l(\Omega')} \le C_0 |\delta|$ for any fixed 
$l$ (by Theorem~\ref{ctythm} part (ii)), it is sufficient to show that 
for some $\alpha_0 >0$ 
\begin{equation*}
	|D\dot{u}|_{C^{k-1,\alpha_0}(\Omega')}
\end{equation*}
is {\em bounded\/} independently of $H$.  Such a bound follows from 
Theorem~$\text{8.32}'$ when applied to the equation in (\ref{dotu}).  
One must check that the coefficients $a_{ij}$ are bounded in 
$C^{k-1,\alpha_0}(\Omega')$ and that $|\dot{u}|_{C^0(\Omega)}$ can also 
be bounded (independently of $H$).   The first bound follows from 
Lemma~\ref{A1} and the bound on $|u|_{C^{k,\alpha_0}(\Omega')}$ given 
in (\ref{nicest}).  The latter bound follows from Theorem~8.16.  
This completes the proof of Lemma~\ref{A}.

\smallskip

From the estimate (\ref{C}) it follows that 
\begin{equation*}
	\ddot{u} = \frac{\partial^2u}{\partial H^2}
\end{equation*}
exists and is well defined in $C^\infty(\Omega)$---satisfying the 
boundary value problem
\begin{equation*}
\begin{cases}
	\sum_{ij} D_i \left( a_{ij} D_j\ddot{u} \right)  = 
		f_2 & \text{on $\Omega$,} \\
	\ddot{v}_{\big|_{\partial\Omega}}  \equiv 0& \text{ }
\end{cases}
\end{equation*}
where $f_2 = -\sum_{i,j} \dot{a}_{ij} D_j\dot{u}$ and 
\begin{equation}\label{dota}
	\dot{a}_{ij} = 
	\sum_k \frac{\partial^2 A}{\partial p_j\partial p_k}(Du)\ D_k\dot{u}.  
\end{equation}

\medskip

Since we have infinitely many derivatives to go (in proving 
Theorem~\ref{regularitythm}), let us assume for $l = 1,2,\ldots,m$, that 
$u^{(l)} = \partial^lu/\partial H^l \in C^\infty(\Omega)$ satisfies 
\begin{equation}\label{induct}
\begin{cases}
	\sum_{i,j} D_i \left( a_{ij} D_ju^{(l)} \right)  = 
		f_{l} & \text{on $\Omega$,} \\
	\text{ }& \text{ } \\
	u^{(l)}_{\big|_{\partial\Omega}}  \equiv 0& \text{ }
\end{cases}
\end{equation}
where $f_1 = 2H$, $f_2 = 2$, and 
\begin{equation}\label{genf}
f_{l+1} = f_{l+1}(x;H) = \dot{f}_l - \sum_{i,j}\dot{a}_{ij} D_j u^{(l)},\quad
2<l \le m-2
\end{equation}
with $\dot{a}_{ij}$ given by (\ref{dota}).

Under these assumptions, we reason as follows:
\begin{theorem}\label{genestthm}
For $|H|,|H_0| \le \tilde{H} < H_{\max}$ and $\Omega''\subset\subset\Omega$, 
$v^{(l)} \equiv u^{(l)} - u_0^{(l)} = 
u^{(l)}(\,\cdot\,;H) - u^{(l)}(\,\cdot\,;H_0)$ satisfies an estimate
\begin{equation}\label{genest}
	|v^{(l)}|_{C^k(\Omega'')} \le C_0 |\delta|
\end{equation}
where $C_0$ is independent of $H$ and $\delta = H-H_0$.
\end{theorem}
\begin{corollary} \label{genind}
There exists $u^{(m+1)} \in C^\infty(\Omega)$ satisfying (\ref{induct}) 
with $m+1$ in place of $l$.
\end{corollary}

Note that Theorem~\ref{regularitythm} clearly follows from 
Corollary~\ref{genind}, and Theorem~\ref{genestthm} is a direct 
generalization of Theorem~\ref{ctythm}.

\medskip

\noindent{\bf Proof of Theorem~\ref{genestthm}.}  According to (\ref{induct}), 
$v^{(m)}$ satisfies
\begin{equation}\label{dill}
\begin{cases}
	\sum_{i,j} D_i \left( a_{ij}^0 D_jv^{(m)} \right)  = 
		f_{l} - f^0_{l} - \sum_{i,j} D_i[ (a_{ij} - a^0_{ij})
			D_ju^{(m)}]& \text{on $\Omega$,} \\
	\text{ }& \text{ } \\
	v^{(m)}_{\big|_{\partial\Omega}}  \equiv 0& \text{ }
\end{cases}
\end{equation}
where $f^0_m = f_m(x;H_0)$ and the other functions 
have been defined above, n.b. (\ref{dotv}).  We take, as before, 
$\Omega'' \subset\subset \Omega' \subset\subset \Omega$, and 
Lemma~$\text{8.32}'$ implies for any $k$ and $\alpha$
\begin{equation*}
	|v^{(m)}|_{C^{k,\alpha}(\Omega'')} \le 
		C_0 (|v^{(m)}|_{C^0(\Omega')} + 
		|f_{m} - f^0_{m}|_{C^{k-1}(\Omega')} + 
		\sum_{i,j} |f_{ij}|_{C^{k-1,\alpha}(\Omega')}  )
\end{equation*}
where $f_{ij} = (a_{ij} - a^0_{ij})D_ju^{(m)}$.

Thus, in general, we have three terms to estimate:
\begin{lemma}\label{vun}
$|f_{ij}|_{C^{k-1,\alpha_0}(\Omega')} \le C_0 |\delta|$, (for some $\alpha_0$).
\end{lemma}
\begin{lemma}\label{tvo}
	$|f_{m} - f^0_{m}|_{C^{k-1}(\Omega')} \le C_0 |\delta|$.
\end{lemma}
\begin{lemma}\label{tve}
	$|v^{(m)}|_{C^0(\Omega')} \le C_0 |\delta|$.
\end{lemma}

Lemmas~\ref{vun} and \ref{tve} follow from the reasoning that gave us 
Lemmas~\ref{A} and \ref{B}, provided we produce a 

\noindent{\bf Proof of Lemma~\ref{tvo}.}  From the definition
\begin{equation*}
	|f_{m} - f^0_{m}|_{C^{k-1}(\Omega')} \le
		|\dot{f}_{m-1} - \dot{f}^0_{m-1}|_{C^{k-1}(\Omega')} +
		\sum_{i,j} |\dot{a}_{ij}D_ju^{(m-1)} - 
			\dot{a}^0_{ij}D_ju_0^{(m-1)} |_{C^{k-1}(\Omega')}.
\end{equation*}
The first term can be handled by induction using the assertion of the 
lemma itself.  If we consider one of the terms in the sum we have
\begin{multline*}
	|\dot{a}_{ij}D_ju^{(m-1)} - 
		\dot{a}^0_{ij}D_ju_0^{(m-1)} |_{C^{k-1}(\Omega')} \\
	\le 
	C_0 ( |\dot{a}_{ij} - \dot{a}^0_{ij} |_{C^{k-1}(\Omega')} 
		|D_ju^{(m-1)}|_{C^{k-1}(\Omega')}  + 
		|\dot{a}^0_{ij}|_{C^{k-1}(\Omega')} 
			|D_jv^{(m-1)}|_{C^{k-1}(\Omega')} ).
\end{multline*}
The first product on the right can be handled by the reasoning in the 
proof of Lemma~\ref{A}.  The second we can estimate by incorporating 
(\ref{genest}) in our induction hypothesis. 

This completes the proof of Lemma~\ref{tvo} and Theorem~\ref{genestthm}.

\smallskip

\noindent{\bf Proof of Corollary~\ref{genind}.}  Since $k$ is arbitrary 
in Theorem~\ref{genestthm}, any sequence $\delta_j \to 0$ provides a 
(sub)sequence of difference quotients $\Delta^\delta u^{(m)} \to w \in 
C^2(\Omega'')$.  Dividing (\ref{dill}) by $\delta$ we also have
\begin{equation*}
	\sum_{i,j} D_i(a^0_{ij} D_j \Delta^\delta u^{(m)}) = 
		\Delta^\delta f_m - \sum_{i,j} D_i(\Delta^\delta a_{ij} 
			D_j u^{(m)} ).
\end{equation*}
On the other hand, it follows inductively from (\ref{genf}) that $f_m$ is 
a linear combination of terms $a_{ij}^{(k)}D_ju^{(l)}$ where 
$a_{ij}^{(k)} = \partial^k a_{ij}/\partial H^k$ and $k,l \ge 1$, $k+l \le m$. 
Consequently, $\dot{f}_m = \lim_{\delta \to 0} \Delta^\delta f_m$ is well 
defined, and the limit $w$ mentioned above satisfies
\begin{equation*}
\begin{cases}
	\sum_{i,j} D_i \left( a_{ij} D_jw \right)  = \dot{f}_{m} - 
	\sum_{i,j} D_i(\dot{a}_{ij} D_ju^{(m)})& \text{on $\Omega$,} \\
	w_{\big|_{\partial\Omega}}  \equiv 0& \text{ }
\end{cases}
\end{equation*}
Since the solutions of this boundary value problem are unique, the limit 
of {\em every} such subsequence must be 
\begin{equation*}
	w= \lim_{\delta\to 0} \Delta^\delta u^{(m)} = u^{(m+1)}.
\end{equation*}
This completes the proof of Corollary~\ref{genind} and, hence, of 
Theorem~\ref{regularitythm}. \Qed

\medskip

We conclude this section with the following observation.
\begin{theorem}[relation of volume and mean curvature]  \label{volthm} Let 
\begin{equation}\label{volume}
	W= W(H) \equiv \int_{\Omega} u(x;H).
\end{equation}
There is a unique value $V_{\max} >0$ such that 
$W: [-H_{\max},H_{\max}] \to [-V_{\max},V_{\max}]$ is a smooth strictly 
decreasing function.
\end{theorem}

\noindent {\bf Proof.}  
\begin{equation*}
	\dot{W} = \int_{\Omega} \dot{u}.
\end{equation*}
Recall from (\ref{dotu}) that $\dot{u}$ satisfies
\begin{equation*}
\begin{cases}
	\sum_{ij} D_i \left( a_{ij} D_j\dot{u} \right)  = 
		2& \text{on $\Omega$,} \\
	\dot{u}_{\big|_{\partial\Omega}}  \equiv 0.& \text{ }
\end{cases}
\end{equation*}
By the maximum principle, any solution of this problem (for any $H$) is 
negative.  \Qed

\section{Convexity for small volumes}

The regularity result of the previous section allows us to linearize the 
problem \astar\  around the zero solution and determine (in terms of 
Poisson's equation) the signs of expressions involving relatively high 
derivatives of $u(x;H)$ in both $x$ and $H$.  If these expressions are 
chosen appropriately as below, we obtain information about the convexity 
of solutions with $|H|$ small.  

In our introductory remarks we were somewhat carefree with the term 
{\em convexity}.  This is essentially justified by the symmetry 
(Corollary~\ref{monotonicity}) of the $H$-graphs under consideration.  
Nevertheless, it will be convenient from now on to distinguish between 
{\em concave\/} functions ($D^2u \le 0$) and {\em convex\/} functions 
($D^2u \ge 0$).  See \cite{MorMul}, Lemma~1.8.1, for equivalent 
definitions.

It will also be convenient for us to detect convexity (or concavity) by 
considering a {\em single number}.  A simple way to do this in two 
dimensions is the following.  Let $v\in C^2(\Omega)$.  We say that 
$v$ is {\em strictly second order convex} (alt. concave) if $D^2v$ is 
positive definite (alt. negative definite) on $\Omega$.  Define an 
auxiliary function $G_v$ on $\Omega$ by 
\begin{equation*}
	G_v = v_{xx}v_{yy} - v_{xy}^2
\end{equation*} 
where we have used the classical ``$x,y$'' notation to denote the second 
partials.  We then have
\begin{lemma}\label{glem}
If $\inf_{\Omega}v < \inf_{\partial\Omega}v$ and $G_v>0$ on $\Omega$, $v$ 
is strictly second order convex.  Similarly, 
$\sup_{\Omega}v > \sup_{\partial\Omega}v$ (and $G_v > 0$) implies $v$ is 
strictly second order concave.
\end{lemma}

\noindent{\bf Proof.}  Since $D^2v(x)$ is a 
real symmetric matrix, there is an orthogonal matrix $M$ and a diagonal 
matrix $\Lambda$ (both of which depend smoothly on $x$) such that 
$MD^2v M^{-1} = \Lambda$.  If the diagonal elements of $\Lambda$ are 
$\lambda_1$ and $\lambda_2$, then $G_v = \text{det}\, D^2v = \lambda_1
\lambda_2$, and the convexity form 
$D^2v e_\theta\cdot e_\theta  = \lambda_1 \xi_1^2 + \lambda_2 
\xi_2^2$ where $\xi = Me_\theta$ and $e_\theta = (\cos\theta,\sin\theta)$.  
Thus, if $G_v > 0$, then neither $\lambda_1$ nor $\lambda_2$ can vanish.  

On the other hand, if $\inf_{\Omega}v < \inf_{\partial\Omega}v$ then there 
is a large lower-hemispherical graph $h\le v$ such that at one or more 
points $h(x_0) = v(x_0)$. Therefore, the $\lambda_i$ must be positive 
(and $v$ strictly second order convex).  \Qed

\medskip

Kawohl \cite{KawWhe} gives essentially the same reasoning as in the above 
proof under the additional assumption that $v$ be subharmonic.

\medskip

The assumptions in Lemma~\ref{glem} are natural for the applications we 
have in mind, but a more general discussion may be found in 
Appendix~\ref{convexityapp}.  At present, our primary objective was to 
justify the following terminology.

\smallskip

We say that $v$ is {\em uniformly} second order convex (alt. concave) if 
$v$ is convex (alt. concave) and for some $\lambda >0$ we have $G_v \ge 
\lambda $ on $\Omega$.  

\medskip

The two ``small volume'' results of this section are obtained by the 
following basic line of reasoning.

If $u=u(x;H)$ is the solution discussed in \S1, then the convexity 
properties of $u$ are the same as those of $u/H$.  To be precise, if 
$w=u/H$ and $G_w >0$, then $G_u >0$.  On the other hand the scaled 
function $w$ is the difference quotient 
\begin{equation}\label{dq}
	\Delta^Hu = \frac{u(x;H) - u(x;0)}{H}
\end{equation}
which, according to the proof of Theorem~\ref{regularitythm} converges in 
$C^2(\bar{\Omega})$ to $\dot{u}$.  Since $G_w$ is a second order operator 
in $x$ we have
\begin{equation*}
	|G_w - G_{\dot{u}}|_{C^0(\bar{\Omega})} \to 0 \quad 
		\text{as $|H|\to 0$.}
\end{equation*}

On the other hand, at $H=0$ we have from (\ref{dotu}) that $\dot{u}$ 
is a solution of the Saint Venant torsion problem
\begin{equation*}
\begin{cases}
	\Delta \dot{u} = 2& \text{on $\Omega$}, \\
	\dot{u}_{\big|_{\partial\Omega}} \equiv 0.
\end{cases}
\end{equation*}
Combining these observations, we have proved
\begin{theorem}[1 convexity]\label{1convex}  Let $\Omega$ be a strictly 
convex domain in the sense that the curvature $\kappa$ of 
$\partial \Omega$ is everywhere positive.  
Consider the problem
\begin{equation}\label{venant}
\begin{cases}
	\Delta v = 2& \text{on $\Omega$}, \\
	v_{\big|_{\partial\Omega}} \equiv 0.
\end{cases}
\end{equation}
If $v$ is uniformly second order convex, then there is some 
$\epsilon >0$ such that 
$u(x;H)$ is strictly second order convex for $0 < H < \epsilon$.  
In particular, if $\Omega$ is an ellipse, ``small bubbles'' are convex.

If on the other hand, $v$ has a point of strict non-convexity, i.e., the 
Gauss curvature of $\text{graph}(v)$ is negative at some point, then 
arbitrarily small bubbles $u$ are likewise non-convex.  This may be observed 
for smooth convex domains whose boundaries converge to a square.
\end{theorem}

Although solutions $\dot{u}$ of (\ref{venant}) are not convex in general, 
they are ``1/2 power convex.''  That is to say, $v=(-\dot{u})^{1/2}$ is 
strictly second order concave.  The crucial step in proving this fact (showing 
that $G_v$ is subharmonic) was carried out by Makar-Limonov \cite{MakSol} 
though the strict second order convexity was actually noted later in 
\cite{KawWhe}.  We note further that this condition is {\em uniform\/}.  In 
order to see this, we introduce another auxiliary function 
$L_v \equiv v_y^2 v_{xx} - 2 v_x v_y v_{xy} +v_x^2 v_{yy}$ which 
essentially measures the convexity of the level curves; see 
Appendix~\ref{convexityapp}.  Given any positive function 
$\phi$ defined on $\Omega$ it is easy to see that 
\begin{equation}\label{F}
	G_\psi = \frac{2\phi G_\phi - L_\phi}{8\phi^2}
\end{equation}
where $\psi = \sqrt{\phi}$.  In our case, 
\begin{equation*}
	G_v = \frac{L_{\dot{u}} - 2\dot{u} G_{\dot{u}} }{8\dot{u}^2}.
\end{equation*}
Now in some closed neighborhood ${\mathcal N}$ of $\partial\Omega$ we may 
assume $|D\dot{u}|\ge \delta >0$ and (consequently) that the curvature of 
the level curves $L_{\dot{u}}/|D\dot{u}|^3 \ge \kappa/2 > 0$.  Taking a 
smaller neighborhood if necessary we may also assume that $|\dot{u}| < 1$ 
and $|\dot{u} G_{\dot{u}}| < \kappa\delta^3/8$.  Thus, on ${\mathcal N} 
\backslash \partial\Omega$ 
\begin{align*}
	G_v &= \frac{1}{8\dot{u}^2} \left( \frac{L_{\dot{u}}}{|D\dot{u}|^3} 
		|D\dot{u}|^3  - 2 \dot{u}G_{\dot{u}} \right)\\
	&\ge \frac{1}{8} \left( \frac{\kappa}{2} 
		\delta^3  - \frac{\kappa\delta^3}{4} \right)\\
	& = \frac{\kappa\delta^3}{32} > 0.
\end{align*}
Finally, on $\Omega' = \Omega\backslash {\mathcal N}$ (by smoothness and 
the observation of Kawohl) there is some $\lambda' >0$ such that 
$G_v\ge\lambda'$.  Letting $\lambda = \min\{\kappa\delta^3/32,\lambda'\}$ 
we have established uniformity.

Our basic line of reasoning now yields
\begin{theorem}\label{lcthm}
Given a strictly convex ($\kappa >0$) domain $\Omega$, there is some 
$\epsilon = \epsilon(\Omega) > 0$ such that $u=u(x;H)$ is $1/2$ concave 
for $-\epsilon < H < 0$.
\end{theorem}

\noindent{\bf Proof.}  Again we scale up.  For the function $w=\sqrt{-u/H}$ 
we have from (\ref{F})
\begin{equation*}
	G_w = \frac{1}{8(\Delta u)^2}(L_{\Delta u} - 2 \Delta u G_{\Delta u})
\end{equation*}
where $\Delta u = \Delta^Hu$ is the difference quotient given in (\ref{dq}).  
Since $\Delta u \to \dot{u}$ in $C^2(\bar{\Omega})$ as $H\to 0$ we conclude 
(essentially from the discussion above) that $L_{\Delta u} - 
2 \Delta u G_{\Delta u} \ge \mu/2 >0$ for 
$|H|$ small enough where.  Thus, $G_w>0$ on $\Omega$ for $|H|$ small.  \Qed

\medskip

Were we able to extend the reasoning of Makar-Limonov to the linear 
problem (\ref{dotu})---and show any degree of strict power convexity---the 
methods of this section would apply to show the convexity of the level 
curves for solutions of \astar.  Another related (and perhaps more 
tractable) approach will be described at the end of the paper.

For now, we concentrate on showing that the ``level curves'' are in 
fact smooth simple closed curves.

\section{Stability}

For this section, let $u = u(x;H)$ be a positive solution $(-H_{\max}<H<0)$ 
of ($\text{\begin{large}$\mathbf\ast\ast$\end{large}}$) on a smooth domain 
in ${\mathbb R}^n$.  We first observe that 
$\text{graph}(u)$ has minimal area among smooth graphs that enclose the same 
volume.  More precisely, if $v\in C^\infty(\Omega)\cap C^0(\bar\Omega)$, 
$v_{\big|_{\partial\Omega}} \equiv \phi$, and 
\begin{equation*}
	V(v) \equiv \int_{\Omega} v = V(u),
\end{equation*}
then
\begin{align}
	A(v) &\equiv \int_{\Omega} \sqrt{1+|Dv|^2} \notag \\
	  & = A(u) + \int_0^1 \frac{d}{dt} \left[ 
		\int_\Omega \sqrt{1+ | (1-t)Du + tDv|^2} \right]\, dt 
		\label{step1} \\
	  & = A(u) + \int_0^1 \left[ \int_\Omega 
		\frac{[(1-t)Du + tDv]\cdot (Dv-Du)}{
		\sqrt{1+ | (1-t)Du + tDv|^2}} \right]\, dt 
		\label{step2} \\
	  & = A(u) + \left[ \int_\Omega 
		\frac{Du \cdot (Dv-Du)}{\sqrt{1+ |Du|^2}} \right] 
		\label{step3}\\
 	  &\qquad\qquad\qquad\qquad	+ \int_0^1 \left[ \int_\Omega 
		\frac{(1+|Dv_*|^2) |Dh|^2 - (Dv_*\cdot Dh)^2}{
		(1+ |Dv_*|^2)^{3/2} } \right]\, dt \notag
\end{align}
where we have expanded the integrand in (\ref{step2}) by Taylor's formula 
at $t=0$; $v_* = (1-t_*)u + t_* v$ for some $t_* \in (0,1)$, and $h= v-u$.  
The numerator in the third term of (\ref{step3}), 
\begin{equation*}
	|Dh|^2 + |Dv_*|^2 |Dh|^2 - (Dv_* \cdot Dh)^2,
\end{equation*}
is nonnegative by Schwarz' inequality, and integrating the second term 
by parts yields
\begin{equation*}
	\int_\Omega \text{div}\,\left(\frac{Du}{\sqrt{1+|Du|^2}}\right) 
		(u-v) = 2H\int_\Omega(u-v) = 0.
\end{equation*}
Hence, $A(v) \ge A(u)$.  (In other words, the area integrand in (\ref{step1}) 
is a convex function of $t$ and has a critical point, hence a minimum, at 
$t=0$.)

Thus, $u$ is a {\em stable critical point for $A$ with 
respect to volume preserving variations\/}:
\begin{corollary}  \label{stabcor}
If $v = v(x;\epsilon)$ is a smooth volume preserving variation of $u$, 
(i.e., $v\in C^\infty(\Omega \times I) \cap C^0(\bar\Omega \times I)$ for 
some interval $I = (-\epsilon, \epsilon)$ and satisfies $V(v) \equiv V(u)$, 
$v\equiv \phi$ on $\partial\Omega$, and $v(x;0)\equiv u$), then 
\begin{equation}\label{fvar}
	\delta A(v) \equiv \frac{d}{d\epsilon} A(v)_{\big|_{\epsilon = 0}} 
			= 0,\quad \text{and}
\end{equation}
\begin{equation}\label{svar}
	\delta^2 A(v) \equiv 
		\frac{d^2}{d\epsilon^2} A(v)_{\big|_{\epsilon = 0}} \ge 0.
\end{equation}
\end{corollary}
\noindent We say that $u$ is a {\em critical point for $A$ with respect to 
volume preserving variations\/} if (\ref{fvar}) holds and a 
{\em (semi)stable critical point with respect to volume preserving 
variations\/} if (\ref{fvar}) and (\ref{svar})
hold.\footnote{Other 
authors have typically considered ``parametric'' variations, but since 
any smooth variation of a graph is locally ``non-parametric'' for small 
$\epsilon$, we lose no generality for the surfaces under consideration.}

We define an alternative functional
\begin{equation*}
	J(v) = A(v) + 2H V(v).
\end{equation*}
It is easy to see that the assertion of statement (\ref{fvar}) in 
Corollary~\ref{stabcor} is equivalent to the following condition.

\smallskip

\begin{equation}
\text{$\delta J(v) = 0$ for {\em any} (not necessarily volume 
	preserving) variation $v$.} 	\tag*{$\text{(\ref{fvar})}'$} 
\end{equation}

\medskip

Following the advice of Bolza \cite{BolVor} and Barbosa and do~Carmo 
\cite{BarSta} we note that a similar equivalence does not hold in general 
for the second variation.  To be precise, if ${\mathcal D}$ is a domain 
in a {\em parametric} surface of constant mean curvature, then we have 
\begin{proposition}[\cite{BarSta}; see also \cite{WenThe}]\label{stablem}
${\mathcal D}$ is stable with respect to volume preserving variations 
if and only if $\delta^2J(\vec{v}) \ge 0$ for all smooth compactly supported 
parametric variations $\vec{v}$ satisfying $\delta V(\vec{v}) = 0$.
\end{proposition}

The reasoning of Cabr\'e and Chanillo in the next section, however, 
essentially requires such an equivalence to hold for non-parametric 
solutions.
\begin{definition}
$u = u(x;H)$ is said to be {\em overstable} (or more accurately 
oversemistable) if $\delta^2 J(v) \ge 0$ for {\em all} compactly 
supported variations $v$.
\end{definition}

\begin{theorem}\label{stability_thm}  Let $u_0 = u(x;H_0)$ where 
$|H_0| < H_{\max}$.  Then $u_0$ is overstable.
\end{theorem}

\begin{remark} Several proofs may be given of Theorem~\ref{stability_thm}.  
Probably the simplest---pointed out to me by C.~Chan 
and H.~Wente---is obtained by repeating the calculation leading to 
Corollary~\ref{stabcor} with $J$ in place of $A$ and an arbitrary 
variation $v$ in place of the volume preserving one.  One then 
has $J(v) \ge J(u)$.  The result also follows---as pointed out by 
R.~Schoen---from the discussion in \cite{FisStr} by noting that 
$g = N_3$ (the 
vertical component of the normal to $\text{\rm graph}(u)$) is a 
positive solution to the equation $\Delta g + \| B \|^2 g = 0$ where 
$\Delta$ denotes the intrinsic Laplacian on $\text{\rm graph}(u)$ and 
$\| B \|^2$ the sum of the squares of the principal curvatures of the 
graph.  The proof presented below demonstrates that $u(x;H)$ provides, 
in some sense, the flow which optimally changes volume.  More precisely, 
given any variation $v$, there is a variation $w$ consisting only of 
members of $\{u(x;H) \}$ such that (\ref{DD}) holds.
\end{remark}

\noindent{\bf Proof or Theorem~\ref{stability_thm}.}  Let $v = v(\,\cdot\,;\epsilon)$ be a 
smooth compactly supported variation.

Since $V(v) \to V(u)$ as $\epsilon \to 0$, we may define 
\begin{equation*}\label{vol2}
	H(\epsilon) = W^{-1}(V(v))
\end{equation*}
where $W$ is given by (\ref{volume}) in Theorem~\ref{volthm}.  

Setting $w = u(\,\cdot\,;H(\epsilon))$ we obtain another 
variation.  Since $V(v) = V(w)$ and $A(w) \le A(v)$ we see that 
$J(w) \le J(v)$ with equality at $\epsilon = 0$.  
Therefore, 
\begin{equation}\label{DD}
	\delta^2 J(w) \le \delta^2 J(v).
\end{equation}

On the other hand, we can compute $\delta^2J(w)$ explicitly.
\begin{align*}
	\frac{d}{d\epsilon} J(w) &= 
		\frac{d}{d\epsilon} \int_\Omega \bigg[ 
		\sqrt{1+|Du(x;H(\epsilon))|^2} + 2H u(x;H(\epsilon)) \bigg] \\
	&= H'(\epsilon)\int_\Omega [2H -2H(\epsilon)] \dot{u}(x;H(\epsilon))
\end{align*}
where $\dot{u} = \partial u/\partial H$ (which is well defined by 
Theorem~\ref{regularitythm}) and we have integrated by parts.  Differentiating 
again and setting $\epsilon = 0$, 
\begin{equation*}
	\delta^2 J(w) = -2H'(0)^2 \int_\Omega \dot{u}(x;H).
\end{equation*}
Since $\dot{u}$ is a solution of (\ref{dotu}) (see the proof of 
Theorem~\ref{volthm}), the integral on the right is strictly negative, and 
$\delta^2 J(v) \ge \delta^2 J(w) \ge 0$.  \Qed

\begin{remark}
The condition of overstability has been considered by various authors 
including Gulliver \cite{GulReg}, Mori \cite{MorSta}, and Ruchert 
\cite{RucEin}.  Ruchert obtains the condition 
\begin{equation}\label{ruc}
	\int_{\mathcal G} \frac{1}{2} \| B\|^2 <2\pi
\end{equation}
for overstability where ${\mathcal G} = \text{graph}(u)$ and $\| B\|^2$ is 
the sum of the squares of the principal curvatures.  Finn has pointed out 
that this condition in inadequate to show the overstability of $H$-graphs 
as follows.
\begin{example}\label{nodex}
Rewriting the integral in (\ref{ruc}) we have 
\begin{align*}
	\int_{\mathcal G} \frac{1}{2} \| B\|^2 &= 
		2\int_{\mathcal G} H^2 - \int_{\mathcal G} K \\
	&\ge H^2 A(u).
\end{align*}

Let $N_r$ be a nodoid (i.e., an inflectionless, rotationally symmetric 
surface of non-zero constant mean curvature) with $H = 1$ and maximum 
distance from its axis of rotation $r$.  If we assume the axis to be the 
$z$-axis and consider ${\mathcal G}_r = \{ (x,y,z)\in N_r: x\ge r-1/8 \}$ 
we obtain graphs with $H^2A_r =  A_r \to +\infty$ as $r\to \infty$.  It 
is clear that (\ref{ruc}) fails for these graphs.  Note: ${\mathcal G}_r$ 
approximates a portion of a circular torus with axis the $z$-axis and 
dimensions $(r-1/2) \times 1/2$.
\end{example}
\end{remark}

\begin{remark}  It should also be noted that while $u(x;H)$ minimizes 
area among graphs that enclose the same volume, it has not been proved 
that $u(x;H)$ is the classical Douglas-Rado-Wente \cite{WenGen} solution 
of the volume constrained Plateau problem.
\end{remark}

Before we proceed, let us recall the formulation of ``overstability'' 
in terms of eigenvalues.  An elementary computation taking $v = u + 
\epsilon \phi$ gives
\begin{equation*}
	\delta^2J(v) = \delta^2A(v) = \langle -{\mathcal L}\phi,\phi \rangle
\end{equation*}
where ${\mathcal L}$ is the linearization of ${\mathcal M}$ at $u$ given in 
(\ref{linearize1}) and the inner product is taken in $L^2$.  Thus, 
the first eigenvalue of ${\mathcal L}$ on $\Omega$, 
$\lambda_1({\mathcal L},\Omega)\ge 0$.  
From the variational characterization of eigenvalues ($\lambda_1 = 
\inf _{|\phi|_{L^2}=1}\langle -{\mathcal L}\phi,\phi \rangle$) and 
the regularity 
of eigenfunctions the following corollary follows at once.
\begin{corollary}\label{e-value}
If $\phi \ne \Omega' \subset\subset \Omega$, then 
$\lambda_1({\mathcal L}, \Omega') > \lambda_1({\mathcal L}, \Omega) \ge 0$.
\end{corollary}

\section{Uniqueness of critical points}\label{cabch}

Here we apply arguments of Cabr\'e and Chanillo \cite{CabSta} to show 
the two theorems stated below.  The reasoning applies to strictly convex 
($\kappa >0$) domains $\Omega$ in ${\mathbb R}^2$.
\begin{theorem}\label{cc1}
If $0 < |H| < H_{\max}$, then for each direction $e_\theta = (\cos\theta,
\sin\theta)$ the solution $u=u(x;H)$ satisfies
\begin{equation}
\begin{cases}
	\text{\rm (i)}\ N_\theta \equiv \{x\in \bar{\Omega}: u_{e_\theta} = 
		Du(x) \cdot e_\theta = 0 \}& \text{is a smooth embedded 
		curve in $\bar{\Omega}$.}  \\
	\text{\rm (ii)}\ M_\theta \equiv \{x\in N_\theta : Du_{e_\theta}
		= 0 \} = \phi.& \text{   }
\end{cases}
\end{equation}
\end{theorem}
\begin{theorem}\label{cc2}
If $u \in C^2(\bar{\Omega})$ is any positive function satisfying 
{\rm (i)}, {\rm (ii)}, and 
\begin{equation*}
\begin{cases}
	u_{\big|_{\partial\Omega}} = 0& \text{ } \\
	|Du_{\big|_{\partial\Omega}}| > 0,& \text{ } 
\end{cases}
\end{equation*}
then $u$ has a unique critical point in $\Omega$.
\end{theorem}

Detailed proofs of both theorems may be found in \cite{CabSta}, but 
for completeness and to give a more detailed exposition of certain points 
we include an outline of the reasoning.

\noindent{\bf Proof of Theorem~\ref{cc1}.}  Notice that $N_\theta$ is a 
smooth embedded curve locally near any point $x_0 \in N_\theta \backslash 
M_\theta$.  Because of this, (i) essentially follows from (ii).  In order 
to verify (ii) we consider two cases.

If $x_0\in M_\theta \cap \partial\Omega$, then $x_0$ must be one of the 
two points $p_1,p_2$ on $\partial\Omega$ where in inward normal $n$ to 
$\partial\Omega$ is orthogonal to $e_\theta$:  $n\cdot e_\theta = 0$.  
(In any case, these two points are in $N_\theta$.)  Calculating the 
normal curvature $\kappa_\sigma$ of ${\mathcal G} = \text{graph}(u)$ 
with respect to $N = (Du,-1)/\sqrt{1+|Du|^2}$ along $\partial\Omega$ 
we have 
\begin{equation}\label{sect}
	\kappa_\sigma = \kappa n\cdot N = 
		\frac{\kappa n \cdot Du}{\sqrt{1+|Du|^2}} = 
		\frac{\kappa}{\sqrt{1+|Du|^2}} \frac{\partial u}{\partial n} 
	> 0.
\end{equation}
On the other hand, since $x_0\in N_\theta$ an alternative expression is 
given by 
\begin{equation*}
	\kappa_\sigma = \frac{\partial^2u}{\partial e_\theta^2}(0,0,1) \cdot N 
		= -\frac{u_{e_\theta e_\theta}}{\sqrt{1+|Du|^2}}.
\end{equation*}
Equating the two expressions we find
\begin{equation}\label{curine}
	Du_{e_\theta}\cdot e_\theta = u_{e_\theta e_\theta} = 
		-\kappa \frac{\partial u}{\partial n} <0.
\end{equation}
Evidently, $Du_{e_\theta}\ne 0$ and the first case is complete.

We have only used the equation in $\astar$ when we asserted, by 
the Hopf boundary point lemma, the inequality in (\ref{sect}).  The 
statement (\ref{curine}) also implies that $N_\theta$ is {\em transverse} 
to $\partial\Omega$ at $x_0 = p_1,p_2$.  

The second possibility is that $x_0 \in M_\theta \cap \Omega$.  This 
assumption may be slightly refined as follows.  From the first case 
there is a closed neighborhood ${\mathcal N}$ of $\partial\Omega$ for which 
${\mathcal N} \cap N_\theta$ consists precisely of connected portions 
$\Gamma_1$ and $\Gamma_2$ of the smooth curves (in $N_\theta$) near 
$p_1$ and $p_2$.  We can also assume that ${\mathcal N}\backslash(\Gamma_1 
\cup \Gamma_2)$ has exactly two connected components $C_+$ and $C_-$ with 
$u_{e_\theta} > 0$ on $C_+$ and $u_{e_\theta} < 0$ on $C_-$.  

Accordingly, we assume $x_0 \in M_\theta \cap(\Omega\backslash{\mathcal N})$.  
According to Hartman \cite[pg. 381 (iv)]{HarEll}, since 
$u_{e_\theta}(x_0) = 0$, $Du_{e_\theta}(x_0) = 0$, and 
${\mathcal L}u_{e_\theta} = 0$, 
a small disk $B$ centered at $x_0$ consists of $4k$ disjoint regions 
in $\Omega\backslash N_\theta$ ($k\ge 1$) along with $4k$ arcs in $N_\theta$ 
connecting $x_0$ to $\partial B$ which are smooth and (with the exception 
of $x_0$) disjoint.  We may furthermore order the regions consecutively 
(say clockwise) so that the first three are $R_+$, $R_-$ and $R_+'$ 
with $u_{e_\theta}$ alternating in sign on the regions as indicated.  
Each region must belong to a connected component of $\Omega\backslash 
N_\theta$.  If every such component $C$ extends to $\partial\Omega$, then 
each such $C$ must be path connected to $C_+$ or to $C_-$.  Assuming this, 
and connecting $R_+$ and $R_+'$ to $C_+$ by paths in 
$\{ x: u_{e_\theta}(x)>0\}$, we see that it is impossible to connect $R_-$ to 
$C_-$ by a path in $\{x: u_{e_\theta}(x) < 0 \}$.  Consequently, some 
component $C = \Omega'$ of 
$\Omega\backslash N_\theta$ is compactly contained in $\Omega$.  Furthermore, 
since $\partial\Omega' \subset N_\theta$, $u_{e_\theta}$ is a nontrivial 
eigenfunction for ${\mathcal L}$:
\begin{equation*}
\begin{cases}
	{\mathcal L}u_{e_\theta} = 0& \text{on $\Omega'$,} \\
	u_{\big|_{\partial\Omega'}} \equiv 0.& \text{ }
\end{cases}
\end{equation*}
This implies that $0 \ge \lambda_1({\mathcal L},\Omega')$ and contradicts 
Corollary~\ref{e-value}.  The contradiction establishes Theorem~\ref{cc1}.  
\Qed

\noindent{\bf Proof of Theorem~\ref{cc2}.}  The basic assertion in this 
proof is that there is a natural flow $\xi:\bar{\Omega} \times {\mathbb R} 
\to \bar{\Omega}$ that ``rotates'' the nodal sets $N_\theta$.  That is, we 
primarily want $\xi$ to satisfy the condition 
\begin{equation}\label{rotat}
	\xi(N_\theta,\tau) = N_{\theta + \tau}.
\end{equation}
This flow also fixes the set of critical points $K \equiv \{x\in \Omega:
Du(x) = 0\} = \cap N_\theta$.  In terms of an autonomous system of ode's 
\begin{equation*}
\begin{cases}
	\dot{\xi} = \vec{F}(\xi),& \text{ } \\
	\xi(x,0) = x,& \text{ } 
\end{cases}
\end{equation*}
the condition on the critical points becomes
\begin{equation*}\label{krit}
	\vec{F}(x) \equiv 0, \quad x\in K.
\end{equation*}

\medskip

At points $x$ away from the critical set $K$, (\ref{rotat}) imposes a 
useful necessary condition on $\vec{F}$ as follows.  Let $\bar{\theta} = 
\bar{\theta}(x)$ be defined by 
\begin{equation*}
	e_{\bar{\theta}} = (\cos\bar{\theta},\sin\bar{\theta}) = 
		(u_y, -u_x)/|Du|.
\end{equation*}
Notice that $x\in N_{\bar{\theta}} \backslash N_\theta$ for $\theta \ne 
\bar{\theta}$.  Consequently, as long as $\xi(x,\tau) \notin K$ we conclude 
from (\ref{rotat}) that 
\begin{equation*}
	\bar{\theta}(\xi(x,\tau)) = \bar{\theta}(x) + \tau.
\end{equation*}
Differentiating with respect to $\tau$,
\begin{equation}\label{ropri}
	D\bar{\theta}\cdot \vec{F} \equiv 1.
\end{equation}

On the other hand, we can compute $D\bar{\theta}$ explicitly:
\begin{equation*}
	D\bar{\theta} = \frac{1}{|Du|^2} D^2u\cdot (-u_y,u_x) = 
		-\frac{1}{|Du|} D^2u\cdot e_{\bar{\theta}}.
\end{equation*}
Looking then at (\ref{ropri}), there is an obvious choice for $\vec{F}$:
\begin{equation*}
	\vec{F}_0 = \frac{D\bar{\theta}}{|D\bar{\theta}|^2} = 
		-|Du|
	\frac{D^2u\cdot e_{\bar{\theta}}}{|D^2u\cdot e_{\bar{\theta}}|^2}.
\end{equation*}
Note that the condition (ii) $M_\theta = \phi$ implies that  
$|D^2u\cdot e_{\bar{\theta}}| = |Du_{e_\theta}| \ge \lambda >0$ uniformly 
for points $x\in K$.  By continuity, a similar bound holds in a neighborhood 
of $K$, and one obtains the estimate $|\vec{F}_0| \le C |Du|$ for some 
constant $C$.  From this it follows that $\vec{F}_0$ extends to a Lipschitz 
vector field on $\bar{\Omega}$ that vanishes on $K$.  

Unfortunately, there is no reason be believe that the resulting flow 
leaves $\bar{\Omega}$ invariant, or equivalently that $\vec{F}_0$ is 
proportional to $e_\theta$ on $\partial\Omega$.  There are many other 
choices for $\vec{F}$ however.  In fact, if $\vec{V}$ is (almost) {\em any} 
vector field, then
\begin{equation}\label{vec}
	\vec{F} = \frac{\vec{V}}{D\bar{\theta}\cdot \vec{V}}
\end{equation}
will satisfy (\ref{ropri}) and imply the main condition (\ref{rotat}).  
$\vec{F}_0$ is obtained by (extending to $K$) the particular choice 
$\vec{V}_0 = D\bar{\theta}$.  Another choice, at least near $\partial\Omega$, 
is given by 
\begin{equation*}
	\vec{V}_1 = e_{\bar{\theta}}
\end{equation*}
which---if it can be extended---will ensure invariance of the domain.  
Note first of all that the formula (\ref{vec}) is valid (i.e., finite valued) 
near $\partial\Omega$.  In fact, 
\begin{equation}\label{check}
	D\bar{\theta}\cdot e_{\bar{\theta}} = 
		-\frac{1}{|D\bar{\theta}|} Du_{e_{\bar{\theta}}} \cdot 
			e_{\bar{\theta}} > 0 \quad \text{by (\ref{curine}).}
\end{equation}
By taking a partition of unity:
\begin{equation*}
\begin{cases}
	\phi_0,\phi_1& \text{smooth, nonnegative on 
		$\Omega$, $\sum \phi_j \equiv 1$,} \\
	\phi_0 \equiv 1& \text{on a large convex domain $\Omega'\subset\subset 
		\Omega$, $K\subset \Omega'$,}\\
	\phi_1 \equiv 1& \text{on a small neighborhood ${\mathcal N} 
		\subset\subset \bar{\Omega}\backslash\Omega'$, 
		$\partial\Omega\subset{\mathcal N}$,} 
\end{cases}
\end{equation*}
and considering $\vec{V} = \sum\phi_j\vec{V}_j$ we get the advantages of 
both $\vec{V}_0$ and $\vec{V}_1$.  (Notice that the sign of the inner 
product in (\ref{check}) agrees with $D\bar{\theta}\cdot D\bar{\theta} > 0$.)  

The resulting flow $\xi$ satisfies all the requirements outlined at the 
beginning of the proof, and $N_\theta\ni x \mapsto \xi(x,\pi) \in 
N_{\theta+\pi} = N_\theta$ is a homeomorphism that reverses the endpoints 
$p_1$ and $p_2\in\partial\Omega$ of $N_\theta$.  Such a map has a unique 
fixed point, and this establishes Theorem~\ref{cc2}.  \Qed

Under the conditions established by Theorems~\ref{cc1} and ~\ref{cc2}---in 
particular, that $Du(x;H)$ vanishes at a unique point 
$x_0 \in \Omega$---it follows from Theorem~\ref{s2} 
Appendix~\ref{convexityapp} that the 
convexity of the level curves is equivalent to the condition $H L_u \ge 0$. 
This observation along with the nondegeneracy of the critical point at 
$x_0=x_0(H)$ (see Lemma~\ref{ex}) suggests the following strategy for 
proving the convexity of the level curves.

Let 
\begin{align*}
	H_0&=\max\{H:HL_u \ge 0\}\\
	   &=\max\{H:L_u \ge 0\}.
\end{align*}
Note that $H_0 > 0$ by Theorem~\ref{lcthm}.  If $H_0 < H_{\max}$, 
then $u_0=u(\ \cdot\ ;H_0)$ satisfies for some neighborhood
$\Omega'$ of $x_0=x_0(H_0)$
\begin{equation*}
\begin{cases}
	L_{u_0} > 0& \text{on $\Omega'\backslash \{x_0\}$}.\\
	L_{u_0} > 0& \text{on $\partial\Omega$}.\\
	L_{u_0} = 0& \text{at some point $x_1 \in \Omega \backslash 
		\Omega'$}.
\end{cases}
\end{equation*}
Under these conditions it is natural to try to show 
$L_{u_0}/ |Du_0|^3$ or (more likely) $L_{u_0} /
|Du_0|^2$ is a supersolution in $\Omega \backslash \{x_0 \}$ of some 
homogeneous elliptic equation.  Convexity of the level curves for all 
$0<|H|<H_{\max}$ would follow in either case.

\appendix

\section{Short time existence}\label{appthm1}
We now return to the proof of Theorem~\ref{solvabilitythm}.  We must 
address the extremal case $H=H_{\max}$.  

First of all note that $H_{\max} < \infty$.  In fact, if $u=u(x;H)$ is 
any solution to $\astar$ then by integrating the equation we have
\begin{equation*}
	2H|\Omega| = 
	\int_\Omega \text{div}\,\left(\frac{Du}{\sqrt{1+|Du|^2}}\right)  =
	\int_{\partial\Omega}\nu\cdot N
\end{equation*}
where $\nu$ is the outward pointing unit normal to $\Omega$ and $N$ is the 
normal to $\text{graph}(u)$.  The integral on the right is clearly bounded 
in absolute value by $|\partial\Omega|$.

From Corollary~\ref{monotonicity} it is natural to define
\begin{equation}\label{limdef}
u(x;\pm H_{\max}) \equiv \pm \lim_{H\nearrow H_{\max}} u(x,H).
\end{equation}
Serrin \cite{SerSur} gives a more general bound for $|u|_{C^0(\Omega)}$ 
than that described in \S\ref{s1}.  To be precise he shows 
$|u|_{C^0(\Omega)} = \sup_\Omega |u| \le 1/(\sigma |H|)$ 
for $|H| \ne 0$ and the alternative bound 
\begin{equation*}
	|u|_{C^0(\Omega)} \le 	1/(\sigma |H|) - \sqrt{1/(\sigma |H|)^2-a^2}
\end{equation*}
if $\Omega$ happens to be contained in a disk of radius 
$a \le 1/(\sigma |H|)$.  Since any domain is contained in a disk of radius 
$C(n){\text{diam}}(\Omega)$ for some constant $C(n) < 1$, we have the absolute 
bound
\begin{equation}\label{serbd}
	|u|_{C^0(\Omega)} \le C(n){\text{diam}}(\Omega).
\end{equation}
From (\ref{serbd}) and the monotonicity it is clear that (\ref{limdef}) 
gives a well defined finite pointwise limit satisfying the boundary 
condition of $\astar$.  
In order to show that the equation is satisfied we restrict to a smooth 
domain $\Omega' \subset\subset\Omega$ (n.b. Problem~6.9) and apply 
Corollary~16.7 \cite[pg. 407]{GilEll}.  In our case, 
$\{u(\,\cdot\,;H)\}_{|H| < H_{\max}}$ is bounded in $C^k(\bar{\Omega}')$ 
for any $k$, and by Lemma~6.36 there is a subsequence $u(\,\cdot\,;H_j)$ 
with $H_j\nearrow H_{\max}$ converging to $u(\,\cdot\,;H_{\max})$ in 
$C^k(\bar{\Omega}')$.  Passing to a limit in the equation, we see that 
$\astar_{H_{\max}}$ is satisfied by the limit function 
$u(\,\cdot\,;H_{\max})$.  

If the gradient blow-up condition 
\begin{equation*}
	\sup_{x\in\Omega} |Du(x;H_{\text{max}})| = \infty
\end{equation*}
were to fail, then additional regularity follows from Lemma~6.18---see 
Remark~\ref{rem2} equation (\ref{linform})---and we have 
$u(\,\cdot\,;H_{\max}) \in C^{2,\alpha}(\bar{\Omega})$.  Thus, we 
arrive at a contradiction (of the definition of $H_{\max}$) from 

\begin{theorem}[short time existence] \label{shorttime}
If $u_0 = u(x;H_0)\in C^2(\bar{\Omega})$ solves $\astar_{H_0}$, 
then there is some $\delta>0$ such that $\astar_H$ is solvable 
in $C^2(\bar{\Omega})$ for $|H-H_0| < \delta$.
\end{theorem}

\noindent{\bf Proof.}  Recall that, using Corollary~\ref{monotonicity} and 
the comparison principle, it is enough to assume $H_0 >0$ and find an 
{\em apriori} 
gradient bound for solutions $u(x;H)$ for $H$ in some interval 
$H_0 < H < H_0 + \delta$.  
A common way to obtain such a bound is to produce a barrier, i.e., we want 
to find a fixed value $H > H_0$ and a fixed 
subsolution $w \in C^2(\Omega) \cap C^1(\bar\Omega)$ for the problem 
$\astar_H$.  We produce such a subsolution as follows.

Consider the linearization of ${\mathcal M}$ at $u_0$
\begin{equation*} 
	{\mathcal L}w = \sum_{ij}\left[ A_{ij}(Du_0) D_iD_j w + 
	\sum_k \frac{\partial^2 A_i}{\partial p_j\partial p_k}(Du_0)
		D_iD_ju_0 D_kw 
		\right].
\end{equation*}
We define $w=w(x;H)$ as the solution to the linear boundary value problem 
\begin{equation}
\begin{cases}\label{linprob}
	{\mathcal L}w &= 2H\quad \text{on $\Omega$,} \\
	w_{\big|_{\partial\Omega}} &\equiv \phi.
\end{cases}
\end{equation}
Note that 
\begin{equation}\label{Hzero}
w(x;H_0) \equiv u_0.  
\end{equation}
Furthermore, $w \in C^\infty(\Omega \times {\mathbb R})$.  See 
the proof of Theorem~\ref{regularitythm}.  In particular, $\dot{w} = 
\partial w/\partial H$ satisfies
\begin{equation}\label{wdot}
	\sum_{ij}\left[ A_{ij}(Du_0) D_iD_j \dot{w} + 
	\sum_k \frac{\partial^2 A_i}{\partial p_j\partial p_k}(Du_0)
		D_iD_ju_0 D_k\dot{w} 
		\right] = 2.
\end{equation}

It remains to show that for some $H > H_0$, $w(x;H)$ is a barrier, i.e., 
that ${\mathcal M}w \ge 2H > 2H_0$.  We see immediately that such a 
constant $H$ exists by differentiating ${\mathcal M}w$ with respect to 
$H$ and setting $H=H_0$; using (\ref{Hzero}) and (\ref{wdot}), the value 
is $2$.  \Qed

\medskip

Applying Theorem~\ref{shorttime} to the zero solution $u(x;0)$ we 
obtain the nondegeneracy $H_{\max} > 0$, and the 
gradient blow-up condition follows as outlined above.  This completes 
the proof of Theorem~\ref{solvabilitythm}.

\begin{remark}
We have assumed no convexity in Theorem~\ref{shorttime}, though the 
Wiener condition is required for the linear problem (\ref{linprob}) 
to be solvable with regular boundary values.  This, however, is 
accomplished for us by a simple regularity assumption on 
$\partial\Omega$ (see the discussion preceding Theorem~6.13).  The 
regularity assertion $w\in C^\infty(\Omega \times {\mathbb R})$ is 
completely analogous to Theorem~\ref{regularitythm} and follows from 
the Schauder estimates (without any additional work).
\end{remark}

\section{H\"older inequalities}\label{apphol}
\begin{lemma}\label{A1}
	Let $g:\Omega\to{\mathbb R}^n$ and $V:{\mathbb R^n}\to{\mathbb R}$.  
	If $|V|_{C^{k+1}(G)} \le C_1$ where $G=g(\Omega)$ and 
	$|g|_{C^{k,\alpha}(\Omega)} \le C_2$, then 
	\begin{equation*}
		|V\circ g|_{C^{k,\alpha}(\Omega)} \le B 
	\end{equation*}
	where $B = B(C_1,C_2,k,n)$.
\end{lemma}
To prove Lemma~\ref{A1} we use induction starting from $k=0$ and 
\begin{lemma}
If $u,v: \Omega \to {\mathbb R}$, 
then for any multiindex $\beta$ with $|\beta| = m$,
\begin{equation*}
	[D^\beta(uv)]_{C^\alpha(\Omega)} 
		\le 2^{m+1} |u|_{C^{m,\alpha}(\Omega)} 
			    |v|_{C^{m,\alpha}(\Omega)}.
\end{equation*}
\end{lemma}

\section{Hadamard type theorems}\label{convexityapp}

In what follows $v\in C^2(\bar{\Omega})$ is a positive function with 
zero boundary values on the planar $C^2$ domain $\Omega$.

For lack of a proof, we begin with a conjecture.
\begin{conjecture}
\begin{description}
\item[{\rm (i)}] ${\mathcal G}\equiv \text{graph}(v)$ is concave if and only if 
\begin{equation*}
	G_v \equiv v_{xx}v_{yy} - v_{xy}^2 \ge 0 \quad 
		\text{for $x\in \Omega$.} 
\end{equation*}
\item[{\rm (ii)}] $v$ is $\,-\infty$ concave, 
i.e., $\text{graph}(v)$ has convex level sets $\{x\in\Omega:v(x)>c\}$, if 
and only if
\begin{equation*}
	L_v \equiv v_y^2v_{xx} - 2 v_xv_yv_{xy} +v_x^2 v_{yy} \le 0 \quad 
		\text{for $x\in \Omega$.} 
\end{equation*}
\end{description}
\end{conjecture}

\begin{remark}  The Gauss curvature of $\mathcal G$ is given by 
$G_v/(1+|Dv|^2)^{3/2}$.  The curvature of the level curves (at least 
where $|Dv| \ne 0$) is given by $-L_v/|Dv|^3$.
\end{remark}

The reasoning in the proof of Lemma~\ref{glem} proves most of 
statement (i).  Recall that $\lambda_1$ and $\lambda_2$ denoted the 
eigenvalues of $D^2u$.  If $v$ is concave then $\lambda_1,\lambda_2 \le 0$, 
and we clearly have that $G_v \ge 0$.  If the condition 
$\lambda_1,\lambda_2 < 0$ holds throughout $\Omega$, we say that $v$ is 
{\em strictly second order concave\/}.  In such a case it is clear that 
$G_v > 0$, and the converse is also true.  

The outstanding case of statement (i), that $G_v \ge 0$ on $\Omega$ 
implies $v$ is concave, is related to the borderline case of Hadamard's 
Theorem for Ovaloids.  See \cite{HopDif}, pp. 119-122, esp. Remark~1.5 and 
\cite{doCDif}, pg. 387, Remark~3.  We note, however, that none of the 
proofs of Hadamard's Theorem given in the above references apply in 
a straightforward way to yield what we find by the simple reasoning above.  
The essential difficulty is that the image of the Gauss map need not be 
simply connected.  For the same reason, the argument of Chern and 
Lashoff \cite{CheTot} 
(also for compact surfaces) is unlikely to settle the 
conjecture easily.

Statement (ii) is, presumably, even more difficult.  We have only the 
following restricted version which is used in our remarks at the end of 
\S\ref{cabch}. 
\begin{theorem}\label{s2}
If $v$ has a unique critical point, then statement {\rm (ii)} holds.
\end{theorem}
\noindent{\bf Proof.}  From our assumption of a unique critical point, it 
follows that each of the level sets $\Omega_c = \{ x\in \Omega: v(x) >c\}$ 
for $0\le c < \max v$ is bounded by a smooth simple closed curve.  Since 
the inward normal to this curve is given by $n = Dv/|Dv|$, it is easy to check 
that the curvature with respect to $n$ is given by $-L_v/|Dv|^3$.  Thus, we 
need only show a version of Hadamard's theorem for simple closed planar 
curves.  A proof for this may be found in \cite{doCDif}; see Proposition~1 
pg.~397. \Qed

\medskip

Our reasoning at the end of \S~\ref{cabch} also uses the following 
simple observation.
\begin{lemma}\label{ex}
If $x_0\in\Omega$ is the unique nondegenerate critical point for $v$, 
i.e., $D^2v(x_0)<0$, then in some neighborhood $\Omega'$ of $x_0$, 
$L_v <0$.
\end{lemma}

\noindent{\bf Proof.}  From the nondegeneracy of the critical point, 
we may assume $Dv \ne 0$ and $D^2v <0$ in $\Omega'$.  Since $L_v = 
w^TD^2vw$ where $w = (v_y,-v_x) \ne 0$, our assertion follows at once. \Qed

\bibliographystyle{alpha}
\bibliography{bib/convex}

\end{document}